\documentclass[a4paper,14pt,reqno]{article}
\usepackage{amsmath,amsthm,amssymb}

\usepackage{setspace,titlesec}
\titleformat{\section}[hang]
{\large\bfseries}
{\thesection}
{1em}
\nonstopmode \numberwithin{equation}{section}
\newtheorem{definition}{Definition}[section]
\newtheorem{theorem}{Theorem}[section]
 \newtheorem{corollary}{Corollary}[section]

\newtheorem{lemma}{Lemma}[section]
\newtheorem{example}{Example}[section]
\newtheorem{remark}{Remark}[section]

\newtheorem{proposition}{Proposition}[section]
\allowdisplaybreaks

\begin{document}
\thispagestyle{empty}

\begin{center}
{\large\bf System of split variational inequality problems in semi-inner product spaces}

\renewcommand{\thefootnote}

\vspace{.25in}
{\bf K.R. Kazmi$^1${\footnote{$^1$Corresponding author; E-mail: krkazmi@gmail.com (K.R. Kazmi)}}}
{ \bf and  Mohd Furkan}$^2${\footnote{$^2$E-mail: mohdfurkan786@gmail.com (Mohd Furkan)}}

\vspace{0.15in}
{\it
${}^{}$Department of Mathematics, Aligarh Muslim University, Aligarh 202002, India} \\
\end{center}

\vspace{.7cm}
\baselineskip=18pt
\noindent{\bf Abstract: } We introduce a new system of split variational inequality problems which is a natural extension of split variational inequality problem in semi-inner product spaces. We use the retraction technique to propose an iterative algorithm for computing the approximate solution of the system of split variational inequality problems. Further, the convergence analysis of the iterative algorithm is also discussed. Several special cases which can be obtained from the main result are also discussed.\\



\section{Preliminaries}
\noindent We recall  the following concepts and results, which are needed to define the problem and to prove the main result:

\begin{definition}{\rm\cite{lumer}}\label{s1d1} Let $X$ be a vector space over the field $\mathbb{K} = \mathbb{R}$ (or $\mathbb{C}$) of real (or complex) numbers. A functional $[\cdot,\cdot] : X\times X \to \mathbb{K}$ is called a semi-inner product if it satisfies the following conditions:
\begin{enumerate}
\item[{\rm(1)}] $[x+y,z]=[x,z]+[y,z], \forall x,y,z\in X;$
\item[{\rm(2)}] $[\lambda x,y]=\lambda[x,y], \forall \lambda \in \mathbb{K} ~ and~ \forall x,y\in X;$
\item[{\rm(3)}] $[x,x] > 0, ~for~ x \neq 0;$
\item[{\rm(4)}] $|[x,y]|^2 \leq [x,x][y,y], \forall x,y\in X.$
\end{enumerate}
\end{definition}
The pair $(X,[\cdot,\cdot])$ is called a semi-inner product space.
As it is observed in {\rm\cite{lumer}} that $\|x\|=[x,x]^{1/2}, \forall x\in X$, is a norm on $X$. Hence every semi-inner product space is a normed linear space. On the other hand, in a normed linear space, one can generate semi-inner product in infinitely many different ways. Further, it is noted that a Hilbert space $H$ can be made into a semi-inner product space, while a semi-inner product  is an inner product if and only if the norm it induces verifies the parallelogram law.\\

Let $Y$ be a semi-inner product space and let $T:X\to Y$ be an arbitrary operator.

\begin{definition}{\rm\cite{pap}}\label{s1d2} The generalized adjoint operator $T^+$ of an operator $T$ is defined as follows: The domain $D(T^+)$ of $T^+$ consists of those $y\in Y$ for which there exists $z \in X$ such that
\begin{equation}
[Tx,y]_{Y}=[x,z]_{X}\nonumber
\end{equation}
for each $x\in X$ and $z=T^+y$.
\end{definition}

\begin{remark} $T^+$ is an operator from $D(T^+)$ into $X$ with the nonempty domain $D(T^+)$, since $0\in D(T^+)$. Hence $T^+(0)=0$. As it is observed in [3] that if $X$ and $Y$ are Hilbert spaces then the generalized adjoint operator is the usual adjoint operator. In general, $T^+$ is not linear even for $T$ is a bounded linear operator.
\end{remark}

Let $C$ be a nonempty closed and convex subset of a Banach space $E$. Let $E^*$ be the dual space of $E$ and $\langle\cdot,\cdot\rangle$ denote the pairing between $E$ and $E^*$. The normalized duality mapping $J:E\to 2^{E^*}$ is defined by
\begin{equation}
J(x)=\{f\in E^*:\langle x,f\rangle=\|x\|^2, \|f\|=\|x\|\}\nonumber
\end{equation}
for all $x\in E$.  We denote by $j$ the single normalized duality mapping, i.e., $j(x)\in J(x),~x\in E$.

\begin{definition}{\rm\cite{xu}}\label{s1d3} Let $U= \{x\in E:\|x\|=1\}$. A Banach space $E$ is said to be:
 \begin{enumerate}
\item[{\rm(1)}] uniformly convex if, for any $\epsilon \in (0,2]$, there exists $\delta>0$ such that for any $x,y\in U$,
\begin{equation}
\|x-y\|\geq\epsilon ~ implies~ \left\|\frac{x+y}{2}\right\|\leq1-\delta ;\nonumber
\end{equation}

\item[{\rm(2)}] smooth if the limit $\lim_{t\to 0}\frac{\|x+ty\|-\|x\|}{t}$ exists for all $x,y\in U$;

\item[{\rm(3)}] uniformly smooth if the limit is attained uniformly for $x,y\in U$.

\end{enumerate}
\end{definition}

\begin{definition}\label{s1d4} The modulus of smoothness of a Banach space $E$ is defined by
\begin{equation}
\rho(\tau)={\rm sup}\left\{\frac{1}{2}(\|x+y\|+\|x-y\|)-1:x,y\in X, \|x\|=1, \|y\|=\tau\right\},\nonumber
\end{equation}
where $\rho:[0,\infty)\to [0,\infty)$ is a function.
\end{definition}

\begin{remark}
$E$ is uniformly smooth if and only if $\lim_{t\to 0}\frac{\rho(\tau)}{\tau}=0$. If $E$ is smooth then normalized duality mapping $J$ is single-valued and if $E$ is uniformly smooth then $J$ is uniformly norm to norm continuous on bounded subsets of $E$. If $E$ is a Hilbert space then $J=I$, where $I$ is the identity mapping.
\end{remark}

In 1967, Giles {\rm\cite{giles}} proved that if the underlying semi-inner product space $X$ is a uniformly convex smooth Banach space then it is possible to define a semi-inner product uniquely which has the following properties:
\begin{enumerate}
\item[{\rm(i)}] $[x,y]=0$ for some $x, y\in X$ if and only if $y$ is orthogonal to $x$, i.e., if and only if $\|y\| \leq \|y+\lambda x\|$, for all scalars $\lambda$.
\item[{\rm(ii)}] The semi-inner product is continuous, i.e., for each $x,y\in X$, we have ${\rm Re}[y,x+\lambda y]\to {\rm Re}[y,x]$ as $\lambda \to 0$.
\item[{\rm(iii)}] The semi-inner product is with the homogeneity property, i.e.,
\begin{equation}
[x,\lambda y]=|\lambda|[x,y], \forall \lambda\in \mathbb{K} ~ {\rm and}~ \forall x,y\in X.\nonumber
\end{equation}
\item[{\rm(iv)}] Generalized Riesz representation theorem: If $f$ is continuous linear functional on $X$ then there is a unique vector $y\in X$ such that $f(x)=[x,y], \forall x\in X$.
\end{enumerate}

The sequence space $l^p, p > 1$ and the function space $L^p, p > 1$ are uniformly convex smooth Banach spaces. More precisely, $L^p$ is min$\{p,2\}$-uniformly smooth for every $p>1$. So one can define semi-inner product on these spaces uniquely.

\begin{example}{\rm\cite{giles}} The real sequence space $l^p$ for $1<p<\infty$ is a semi-inner product space with the semi-inner product defined by
\begin{equation}
[x,y]=\frac{1}{\|y\|^{p-2}_p}\sum_{i}x_iy_i|y_i|^{p-2}, \forall x,y\in l^p.\nonumber
\end{equation}
\end{example}

\begin{example}{\rm\cite{giles}} The real Banach space $L^p(X,\mu)$ for $1<p<\infty$ is a semi-inner product space with the semi-inner product defined by
\begin{equation}
[f,g]=\frac{1}{\|g\|^{p-2}_p}\int_{X}fx|gx|^{p-1}sgn(gx)d\mu, \forall f,g\in L^p.\nonumber
\end{equation}
\end{example}

Now, we summarize the following properties of the generalized adjoint operator from the results given in {\rm\cite{pap}}.


\begin{proposition}\label{s1p1} Let $X$ and $Y$ be $2$-uniformly convex smooth Banach spaces and let $T:X\to Y$ be a bounded linear operator. Then
\begin{enumerate}
\item[{\rm(i)}] $D(T^+)=Y$;
\item[{\rm(ii)}] $T^+$ is bounded, and it holds that
\begin{equation}
\|T^+y\| \leq \|T\|\|y\|, \forall y\in Y.\nonumber
\end{equation}
\end{enumerate}

\end{proposition}

\begin{definition}{\rm\cite{reich}}\label{s1d5} Let $D$ be a subset of $C$ and $Q_C$ be a mapping of $C$ into $D$. Then $Q_C$ is said to be sunny if
\begin{equation}
Q_C(Q_Cx+t(x-Q_Cx))=Q_Cx,\nonumber
\end{equation}
whenever $Q_Cx+t(x-Q_Cx)\in C$ for $x\in C$ and $t\geq 0$.
\end{definition}

\begin{definition}{\rm\cite{reich}}\label{s1d6} A subset $D$ of $C$ is called a sunny nonexpansive retract of $C$ if there exists a sunny nonexpansive retraction from $C$ into $D$.
\end{definition}

The following result describes a characterization of sunny nonexpansive retractions on a smooth  Banach space.

\begin{proposition}{\rm\cite{reich}}\label{s1p2} Let $E$ be a smooth Banach space and let $C$ be a nonempty subset of $E$. Let $Q_C:E\to C$ be a retraction. Then the following are equivalent:
\begin{enumerate}
\item[{\rm(i)}] $Q_C$ is sunny and nonexpansive;
\item[{\rm(ii)}] $\|Q_Cx-Q_Cy\|^2 \leq \langle x-y, J(Q_Cx-Q_Cy)\rangle, \forall x,y\in E$;
\item[{\rm(iii)}] $\langle x-Q_Cx, J(y-Q_Cx)\rangle \leq 0, \forall x\in E, y\in C$.
\end{enumerate}
\end{proposition}

\begin{lemma}{\rm\cite{xu}}\label{s1l1} Let $p>1$ be a real number and $E$ be a smooth Banach space. Then the following statements are equivalent:
\begin{enumerate}
\item[{\rm(i)}] $E$ is $2$-uniformly smooth;
\item[{\rm(ii)}] There is a constant $c>0$ such that for every $x,y\in E$, the following inequality holds
\begin{equation}
\|x+y\|^2\leq \|x\|^2+2\langle y,J(x) \rangle + c\|y\|^2.\nonumber
\end{equation}
\end{enumerate}
\end{lemma}

\begin{remark}\label{s1r1}

\begin{enumerate}
\item  {\rm\cite{giles,lumer,sahu}}: Every normed linear space is a semi-inner product space. In fact by Hahn Banach theorem, for each $x\in E$ there exists at least one functional $f_x\in E^*$ such that $\langle x, f_x\rangle=\|x\|^2$. Given any such mapping $f$ from $E$ into $E^*$, we can verify that $[y,x]=\langle y, f_x\rangle$ defines a semi-inner product. Hence, we can write the inequality given in Lemma {\ref{s1l1}} as
\begin{equation}
\|x+y\|^2\leq \|x\|^2+2[y,x] + c\|y\|^2, \forall x,y\in E.\nonumber
\end{equation}
The constant $c$ is chosen with best possible minimum value. We call $c$ as the constant of smoothness of $E$.
\item The inequalities given in Proposition {\ref{s1p2}} (ii) \& (iii) can be written as
\begin{enumerate}
\item[{\rm(ii)}] $\|Q_Cx-Q_Cy\|^2 \leq [ x-y, Q_Cx-Q_Cy], \forall x,y\in E$;
\item[{\rm(iii)}] $[ x-Q_Cx, y-Q_Cx] \leq 0, \forall x\in E, y\in C$.
\end{enumerate}
\end{enumerate}
\end{remark}

\begin{example}{\rm\cite{sahu}} The function space $L^p$ is $2$-uniformly smooth for $p\geq 2$ and it is $p$-uniformly smooth for $1<p<2$.
If $2\leq p <\infty$, then we have for all $x,y\in L^p$,
\begin{equation}
\|x+y\|^2\leq \|x\|^2+2[y,x] + (p-1)\|y\|^2,\nonumber
\end{equation}
where $(p-1)$ is the constant of smoothness.
\end{example}

Let $E_1$ and $E_2$ be $2$-uniformly convex, smooth Banach spaces and for each $i\in \{1,2\}$; let $C_i\subset E_i$ be a nonempty, closed and convex set and let $J_1:E_1\to 2^{E_1^*}$ and $J_2:E_2\to 2^{E_2^*}$ be the normalized duality mappings. Let $F,G:C_1\to E_1$ and $f,g:C_2\to E_2$ be nonlinear mappings, and let $A:E_1\to E_2$ be a bounded linear operator. We introduce the following system of split variational inequality problems (in short, SSpVIP):
Find $(x_1,y_1)\in C_1\times C_1$ such that
$$\langle \lambda F y_1+x_1-y_1, J_1(z_1-x_1)\rangle\geq 0,~~\forall z_1\in C_1, $$
and such that $(x_2,y_2)$ with $x_2=Ax_1\in C_2, y_2=Ay_1\in C_2$ solves
$$\langle \gamma f y_2+x_2-y_2, J_2(z_2-x_2)\rangle\geq 0,~~\forall z_2\in C_2; $$
$$\langle \lambda G x_1+y_1-x_1, J_1(z_1-y_1)\rangle\geq 0,~~\forall z_1\in C_1, $$
and such that $(x_2,y_2)$ solves
$$\langle \gamma g x_2+y_2-x_2, J_2(z_2-y_2)\rangle\geq 0,~~\forall z_2\in C_2, $$
for any $\lambda, \gamma > 0$.\\

Above SSpVIP  is equivalent to find $(x_1,y_1)\in C_1\times C_1$ such that
$$[ \lambda F y_1+x_1-y_1, z_1-x_1]\geq 0,~~\forall z_1\in C_1, \eqno(1.1)$$
and such that $(x_2,y_2)$ with $x_2=Ax_1\in C_2, y_2=Ay_1\in C_2$ solves
$$[\gamma f y_2+x_2-y_2, z_2-x_2]\geq 0,~~\forall z_2\in C_2; \eqno(1.2)$$
$$[ \lambda G x_1+y_1-x_1, z_1-y_1]\geq 0,~~\forall z_1\in C_1, \eqno(1.3)$$
and such that $(x_2,y_2)$ solves
$$[ \gamma g x_2+y_2-x_2, z_2-y_2]\geq 0,~~\forall z_2\in C_2, \eqno(1.4)$$
for any $\lambda, \gamma > 0$.\\

\noindent {\bf Some special cases:}\\

{\bf 1.} If we set $E_1=H_1,~E_2=H_2$, where $H_1,~H_2$ are Hilbert spaces, then  SSpVIP (1.1)-(1.4) reduces to the following system of split variational inequality problems (SSpVIP) in Hilbert spaces: Find $(x_1,y_1)\in C_1\times C_1$ such that
$$\langle \lambda F y_1+x_1-y_1, z_1-x_1\rangle\geq 0,~~\forall z_1\in C_1, \eqno(1.5)$$
and such that $(x_2,y_2)$ with $x_2=Ax_1\in C_2, y_2=Ay_1\in C_2$ solves
$$\langle \gamma f y_2+x_2-y_2, z_2-x_2\rangle\geq 0,~~\forall z_2\in C_2; \eqno(1.6)$$
$$\langle \lambda G x_1+y_1-x_1, z_1-y_1\rangle\geq 0,~~\forall z_1\in C_1, \eqno(1.7)$$
and such that $(x_2,y_2)$ solves
$$\langle \gamma g x_2+y_2-x_2, z_2-y_2\rangle\geq 0,~~\forall z_2\in C_2, \eqno(1.8)$$
for any $\lambda, \gamma > 0$.\\

{\bf 2.} If we set $F=G, ~f=g,~ \lambda = \gamma,~ y_1=x_1$, then $y_2=x_2$ and hence SSpVIP (1.1)-(1.4) reduces to the following split variational inequality problem (in short, SpVIP): Find $x_1\in C_1$ such that
$$[ Fx_1, z_1-x_1] \geq 0, ~~\forall z_1\in C_1, \eqno(1.9)$$
and such that $x_2=Ax_1\in C_2$ solves
$$[ fx_2, z_2-x_2] \geq 0, ~~\forall z_2\in C_2, \eqno(1.10)$$

{\bf 3.} In Case 2, if $E_1=H_1,~ E_2=H_2$, then SpVIP (1.9)-(1.10) reduces to the split variational inequality problem considered and studied by Censor et al. {\rm\cite{cen3}}. It is worth mentioning that the SpVIP is quite general and permit split minimization between two spaces so that the image of a minimizer of a given function, under a bounded linear operator, is a minimizer of another function. It includes as a special case, the variational inequality problem, the split zero problem and the split-feasibility problem which have already been studied and used in practice as a model in the intensity-modulated radiation therapy planning, see {\rm\cite{cen2,cen1}}. For a further related work, see {\rm\cite{byr,kaz1,kaz5,kazm1,kazm2,mou}}.\\

Further, it is worth mentioning that so far  the iterative approximations of split variational inequality problem and its generalizations  have been studied  in the setting of Hilbert spaces.  Therefore, a natural
question appears as to whether or not one can  study   these problems in setting of  Banach spaces.\\

In this paper, we  use the retraction technique  to propose and analyze an iterative algorithm for computing the approximate solution of SSpVIP (1.1)-(1.4) in $2$-uniformly convex smooth Banach spaces. Further, convergence analysis of the iterative algorithm is discussed. Several special cases which can be obtained from the main result are also discussed. The problems and the results discussed in this paper are new and different from the existing problems and results in the literature.\\

\section{Iterative Algorithms}
By making use of  Proposition {\ref{s1p2}} , we easily observe that SSpVIP (1.1)-(1.4) can be formulated as follows: Find $(x_1,y_1)\in C_1\times C_1$ with $(x_2,y_2)=(Ax_1,Ay_1)\in C_2\times C_2$ such that
$$x_1=Q_{C_1}(y_1-\lambda Fy_1), \eqno (2.1)$$
$$x_2=Q_{C_2}(y_2-\gamma fy_2), \eqno (2.2)$$
$$y_1=Q_{C_1}(x_1-\lambda Gx_1), \eqno (2.3)$$
$$y_2=Q_{C_2}(x_2-\gamma gx_2), \eqno (2.4)$$
for $\lambda,  \gamma > 0$.\\

Based on above arguments, we propose the following iterative algorithm for approximating a solution to SSpVIP (1.1)-(1.4).\\

Let $\{\alpha_n\}\subseteq (0,1)$ be a sequence such that $\sum\limits^{\infty}_{n=1}=\infty$.\\

\noindent {\bf Iterative Algorithm 2.1.} Given $(x^0_1,y^0_1)\in C_1\times C_1$, compute the iterative sequence $\{(x^n_1,y^n_1)\}$ defined by the iterative schemes:
$$a^n_1=Q_{C_1}(y^n_1-\lambda Fy^n_1), \eqno (2.5)$$
$$a^n_2=Q_{C_2}(y^n_2-\gamma fy^n_2), \eqno (2.6)$$
$$b^n_1=Q_{C_1}(x^n_1-\lambda Gx^n_1), \eqno (2.7)$$
$$b^n_2=Q_{C_2}(x^n_2-\gamma gx^n_2), \eqno (2.8)$$
$$x^{n+1}_1=(1-\alpha^n)x^n_1+\alpha^n\left(a^n_1+\rho A^+(a^n_2-Aa^n_1)\right), \eqno (2.9)$$
$$y^{n+1}_1=(1-\alpha^n)y^n_1+\alpha^n\left(b^n_1+\rho A^+(b^n_2-Ab^n_1)\right), \eqno (2.10)$$
for all $n=0,1,2,....$ and $\lambda,  ~\gamma, ~\rho > 0$, where $A^+$ is  the generalized adjoint operator of $A$, and $x^n_2=Ax^n_1$ and $y^n_2=Ay^n_1$ for all $n$.\\

If we set $E_1=H_1, ~E_2=H_2$, where $H_1,~H_2$ are Hilbert spaces, then  Iterative Algorithm 2.1 reduces to the following iterative algorithm for computing the approximate solution of SSpVIP (1.5)-(1.8):\\

\noindent {\bf Iterative Algorithm 2.2.} Given $(x^0_1,y^0_1)\in C_1\times C_1$, compute the iterative sequence $\{(x^n_1,y^n_1)\}$ defined by the iterative schemes:
$$a^n_1=P_{C_1}(y^n_1-\lambda Fy^n_1), \eqno (2.11)$$
$$a^n_2=P_{C_2}(y^n_2-\gamma fy^n_2), \eqno (2.12)$$
$$b^n_1=P_{C_1}(x^n_1-\lambda Gx^n_1), \eqno (2.13)$$
$$b^n_2=P_{C_2}(x^n_2-\gamma gx^n_2), \eqno (2.14)$$
$$x^{n+1}_1=(1-\alpha^n)x^n_1+\alpha^n\left(a^n_1+\rho A^*(a^n_2-Aa^n_1)\right), \eqno (2.15)$$
$$y^{n+1}_1=(1-\alpha^n)y^n_1+\alpha^n\left(b^n_1+\rho A^*(b^n_2-Ab^n_1)\right), \eqno (2.16)$$
for all $n=0,1,2,....$ and $\lambda, ~ \gamma,~ \rho > 0$, where $A^*$ is  the adjoint operator of $A$ with $\|A^*\|=\|A\|$, and  $P_{C_i}$ is the metric projection of $H_i$ onto $C_i$ for each $i \in \{1,2\}$.\\

If we set $F=G, ~f=g, ~\lambda = \gamma,~ y_1=x_1$, then $y_2=x_2$ and hence Iterative Algorithm 2.1 reduces to the following iterative algorithm for computing the approximate solution of SpVIP (1.9)-(1.10):\\

\noindent {\bf Iterative Algorithm 2.3.} Given $x^0_1\in C_1$, compute the iterative sequence $\{x^n_1\}$ defined by the iterative schemes:
$$a^n_1=Q_{C_1}(x^n_1-\lambda Fx^n_1), $$
$$a^n_2=Q_{C_2}(x^n_2-\lambda fx^n_2),$$
$$x^{n+1}_1=(1-\alpha^n)x^n_1+\alpha^n\left(a^n_1+\rho A^+(a^n_2-Aa^n_1)\right), $$
for all $n=0,1,2,....$ and $\lambda,~ \gamma, ~\rho > 0$.

\section{Main Result}
First, we define the following concepts.

\vspace{.3cm}
\begin{definition}\label{s3d1}
A mapping $F:E_1 \rightarrow E_1$ is said to be
\begin{enumerate}
\item[{\rm(1)}] $\alpha$-{\it strongly monotone} if there exists a constant $\alpha>0$ such that
$$[Fx_1-Fy_1,x_1-y] \geq \alpha\Vert x_1-y_1\Vert^2,~~~ \forall x_1,y_1\in E_1;$$
\item[{\rm(2)}] $\beta$-{\it Lipschitz continuous}, if there exists a constant  $\beta > 0$ such that
$$\|Fx_1-Fy_1\|\leq \beta \|x_1-y_1\|,~~~ \forall x_1,y_1\in E_1.$$
\end{enumerate}
\end{definition}

Now, we prove that the sequence of approximate solutions of SSpVIP (1.1)-(1.4) generated by Iterative Algorithm 2.1 converges strongly to the solution of SSpVIP (1.1)-(1.4).

\begin{theorem}\label{s3t1} For each  $i\in\{1,2\},$ let $C_{i}$ be a nonempty, closed and convex subset of $2$-uniformly convex smooth Banach space $E_{i}$ with constant of smoothness $c_i$. Let $F: C_1\to E_1$ be $\alpha_{1}$-strongly monotone  and $\beta_{1}$-Lipschitz continuous;
  let $G:  C_1\to E_1$ be $\alpha_{2}$-strongly monotone  and $\beta_{2}$-Lipschitz continuous;
let $f: C_2\to E_2$ be $\sigma_{1}$-strongly monotone  and $\eta_{1}$-Lipschitz continuous, and
let $g: C_2\to E_2$ be $\sigma_{2}$-strongly monotone and $\eta_{2}$-Lipschitz continuous. Let $A:E_1 \to E_2$  be bounded linear operator. Suppose $(x_1,y_1) \in C_1 \times C_2$ is a solution to SSpVIP(1.1)-(1.4) then the sequence $\{(x_1^n, y_1^n)\}$ generated by Iterative Algorithm 2.1  converges strongly to $(x_1,y_1)$ provided that  the constant $ \lambda>0$  satisfies the condition:
$$\max_{1\leq i\leq 2}\left\{\frac{\alpha_i-\sqrt{\alpha_i^2-c_1\beta_i^2(1-p_i^2)}}{c_1\beta_i^2}\right\}<\lambda < \min_{1\leq i\leq 2}\left\{\frac{\alpha_i+\sqrt{\alpha_i^2-c_1\beta_i^2(1-p_i^2)}}{c_1\beta_i^2}\right\} \eqno (3.1)$$
$$\alpha_i>\beta_i\sqrt{c_1(1 -p_i^2)};~~p_i=\frac{1-m\theta_{i+2}}{1+m};~~~m=\rho\|A^+\|\|A\|;$$
$$\theta_{i+2}=\sqrt{1-2\gamma\sigma_i+c_2\gamma^2\eta_i^2};~~~ \gamma>0.$$
\end
{theorem}

\noindent {\bf Proof.} Given that $(x_1,y_1)$ is a solution of SSpVIP (1.1)-(1.4), that is, $x_1, y_1$ satisfy the relations (2.1)-(2.4).
Since $F:C_1\to E_1$ is $\alpha_1$-strongly monotone and $\beta_1$-Lipschitz continuous, from Iterative Algorithm 2.1 (2.5) and (2.1), we estimate\\

$$\|a^n_1-x_1\| = \|Q_{C_1}(y^n_1-\lambda Fy^n_1)-Q_{C_1}(y_1-\lambda Fy_1)\| \hspace{.95in}$$
$$\leq  \|y^n_1-y_1-\lambda (Fy^n_1-Fy_1)\| \hspace{1in}$$
$$\hspace{1in} \leq  \left(\|y^n_1-y_1\|^2-2\lambda \left[Fy^n_1-Fy_1,y^n_1-y_1\right]+c\lambda^2\|Fy^n_1-Fy_1\|^2\right)^{\frac{1}{2}}$$
$$\leq \theta_1 \|y^n_1-y_1\|, \hspace{1.75in} ~~~\eqno(3.2)$$
where $\theta_1=(1-2\lambda\alpha_1+c_1\lambda^2\beta^2_1)^{\frac{1}{2}}.$\\

Next, since $G:C_1\to E_1$ is $\alpha_2$-strongly monotone and $\beta_2$-Lipschitz continuous, from Iterative Algorithm 2.1 (2.7) and (2.3), we have
$$\|b^n_1-y_1\|=\|Q_{C_1}(x^n_1-\mu Gx^n_1)-Q_{C_1}(x_1-\mu Gx_1)\| \hspace{1in}$$
$$\leq \theta_2\|x^n_1-x_1\|, \hspace{1.75in} ~~~\eqno(3.3)$$
where $\theta_2=(1-2\lambda\alpha_2+c_1\lambda^2\beta^2_2)^{\frac{1}{2}}.$\\

Again, since $f:C_2\to E_2$ is $\sigma_1$-strongly monotone and $\eta_1$-Lipschitz continuous, from Iterative Algorithm 2.1 (2.6) and (2.2), we have
$$\|a^n_2-x_2\| \leq \theta_3\|y^n_2-y_2\|, ~~~\eqno(3.4)$$
where $\theta_3=(1-2\gamma\sigma_1+c_2\gamma^2\eta^2_1)^{\frac{1}{2}}.$\\

Since $g:C_2\to E_2$ is $\sigma_2$-strongly monotone and $\eta_2$-Lipschitz continuous, from Iterative Algorithm 2.1 (2.8) and (2.4), we have
$$\|b^n_2-y_2\| \leq \theta_4\|x^n_2-x_2\|, ~~~\eqno(3.5)$$
where $\theta_4=(1-2\gamma\sigma_2+c_2\gamma^2\eta^2_2)^{\frac{1}{2}}.$\\

Now, using the fact that $A^+$  is bounded, we have
\begin{eqnarray}
\|x^{n+1}_1-x_1\|  &\leq & (1-\alpha^n)\|x^n_1-x_1\|+\alpha^n\|a^n_1-x_1+\rho A^+(a^n_2-Aa^n_1)\|\nonumber\\
&\leq & (1-\alpha^n)\|x^n_1-x_1\|+\alpha^n\|a^n_1-x_1\|+\alpha^n\rho\|A^+\|\|a^n_2-Aa^n_1\|\nonumber\\
&\leq & (1-\alpha^n)\|x^n_1-x_1\|+\alpha^n\|a^n_1-x_1\|+\alpha^n\rho\|A^+\|\left(\|a^n_2-x_2-Aa^n_1+x_2\|\right)\nonumber\\
&\leq & (1-\alpha^n)\|x^n_1-x_1\|+\alpha^n\|a^n_1-x_1\|+\alpha^n\rho\|A^+\|\left(\|a^n_2-x_2\|+\|A\|\|a^n_1-x_1\|\right)\nonumber\\
&=& (1-\alpha^n)\|x^n_1-x_1\|+\alpha^n\theta_1\|y^n_1-y_1\|+\alpha^n\rho\|A^+\|\left(\theta_3\|y^n_2-y_2\|+\|A\|\theta_1\|y^n_1-y_1\|\right)\nonumber\\
&\leq & (1-\alpha^n)\|x^n_1-x_1\|+\alpha^n\theta_1\|y^n_1-y_1\|+\alpha^n\rho\|A^+\|\|A\|\left(\theta_3\|y^n_1-y_1\|+\theta_1\|y^n_1-y_1\|\right)\nonumber\\
&=& (1-\alpha^n)\|x^n_1-x_1\|+\alpha^n\left(\theta_1+\rho\|A^+\|\|A\|(\theta_1+\theta_3)\right)\|y^n_1-y_1\|. \hspace{1.2in} (3.6) \nonumber
\end{eqnarray}

Similarly, we obtain

$$\|y^{n+1}_1-y_1\| \leq  (1-\alpha^n)\|y^n_1-y_1\|+\alpha^n\left(\theta_2+\rho\|A^+\|\|A\|(\theta_2+\theta_4)\right)\|x^n_1-x_1\|. \hspace{.5in} \eqno(3.7)$$

Now, define the norm $||\cdot||_{\star}$ on $E_1\times E_2$ by
$$||(x,y)||_{\star} =||x||+||y||,~~~(x,y)\in E_1\times E_2.$$
We can easily show that $(E_1\times E_2, ||\cdot||_{\star})$ is a Banach space.

\vspace{.3cm}
By making using of (3.6) and (3.7), we have the following estimate:
$$\|(x_1^{n+1},y_1^{n+1})-(x_1,y_1)\|_* =\|x_1^{n+1}-x_1\|+\|y_1^{n+1}-y_1\|\hspace{2.4in}$$
$$\leq (1-\alpha_n)\left(\|x_1^{n}-x_1\|+\|y_1^{n}-y_1\|\right)\hspace{.4in}$$
$$\hspace{.6in} +\alpha_n\left(\theta_1+\rho\|A^+\|\|A\|(\theta_1+\theta_3)\right)\|y^n_1-y_1\|\|y_1^{n}-y_1\|$$
$$\hspace{.08in}+ \alpha_n\left(\theta_2+\rho\|A^+\|\|A\|(\theta_2+\theta_4)\right)\|x^n_1-x_1\|$$
$$\leq (1-\alpha_n)\left(\|x_1^{n}-x_1\|+\|y_1^{n}-y_1\|\right)\hspace{.4in}$$
$$\hspace{.15in} +\alpha_n \max\{k_1, k_2\}\left(\|x^n_1-x_1\|+\|y_1^{n}-y_1\|\right)$$
$$= \left(1-\alpha_n(1-\theta)\right)\|(x_{n},y_{n})-(x,y)\|_*,\hspace{.27in}\eqno(3.8)$$
where $\theta= \max \{k_1,k_2\};~~k_1=\theta_1+m(\theta_1+\theta_3);~~k_2=\theta_2+m(\theta_2+\theta_4);~~m=\rho\|A^+\|\|A\|.$

\vspace{0.3cm}
Thus, we obtain
 $$\|(x_{n+1},y_{n+1})-(x,y)\|_* < \prod\limits^{n}_{r=1}\left(1-\alpha_r(1-\theta)\right) \|(x_1^{0},y_1^{0})-(x_1,y_1)\|_*. \eqno(3.9)$$

It follows from given condition (3.1) on $\lambda$ that $\theta \in (0,1)$. Since $\sum \limits^{\infty}_{n=1} \alpha_n=\infty$ and $\theta \in (0,1)$, it implies  that
$$\lim\limits_{n\to \infty}\prod\limits^{n}_{r=1}\left(1-\alpha_r(1-\theta)\right)=0.$$

 Thus, it follows from  (3.9) that $\{(x_1^{n+1},y_1^{n+1})\}$  converges strongly to $(x_1,y_1)$ as $n \to \infty$, that is, $x_1^n \to x_1$ and $y_1^n \to y_1$ as $n \to \infty$. Further, it follows from (3.2) and (3.3), respectively, that $a_1^n \to x_1$ and $b_1^n \to y_1$  as $n \to \infty$. Hence, it follows from (3.4) and (3.5), respectively, that $a_2^n\to x_2=Ax_1$ and  $b_2^n \to y_2=Ay_1$  as $n \to \infty$.   This completes the proof.

\vspace{.3cm}
Now, we give the following corollaries which are consequences of Theorem {\ref{s3t1}}.

\vspace{.3cm}
If we set $E_1=H_1,~E_2=H_2$, then Theorem {\ref{s3t1}} reduces to the following result  for the convergence analysis of Iterative Algorithm 2.2 for SSpVIP (1.5)-(1.8).

\begin{corollary}\label{s3t2} For each  $i\in\{1,2\},$ let $C_{i}$ be a nonempty, closed and convex subset of real Hilbert space $H_{i}$. Let $F: C_1\to H_1$ be $\alpha_{1}$-strongly monotone  and $\beta_{1}$-Lipschitz continuous;
  let $G: C_1\to H_1$ be $\alpha_{2}$-strongly monotone  and $\beta_{2}$-Lipschitz continuous;
let $f: C_2\to H_2$ be $\sigma_{1}$-strongly monotone  and $\eta_{1}$-Lipschitz continuous, and
let $g: C_2\to H_2$ be $\sigma_{2}$-strongly monotone and $\eta_{2}$-Lipschitz continuous. Let $A:H_1 \to H_2$  be bounded linear operator. Suppose $(x_1,y_1) \in C_1 \times C_2$ is a solution to SSpVIP(1.5)-(1.8) then the sequence $\{(x_1^n, y_1^n)\}$ generated by Iterative Algorithm 2.2  converges strongly to $(x_1,y_1)$ provided that the constant $ \lambda>0$  satisfies the condition:
$$\max_{1\leq i\leq 2}\left\{\frac{\alpha_i-\sqrt{\alpha_i^2-\beta_i^2(1-p_i^2)}}{\beta_i^2}\right\}<\lambda < \min_{1\leq i\leq 2}\left\{\frac{\alpha_i+\sqrt{\alpha_i^2-\beta_i^2(1-p_i^2)}}{\beta_i^2}\right\}$$
$$\alpha_i>\beta_i\sqrt{(1 -p_i^2)};~~p_i=\frac{1-m\theta_{i+2}}{1+m};~~~m=\rho\|A\|^2;$$
$$\theta_{i+2}=\sqrt{1-2\gamma\sigma_i+\gamma^2\eta_i^2};~~~ \gamma>0.$$
\end
{corollary}

\vspace{.2in}
If we set $F=G, ~f=g,~ \lambda = \gamma, ~y_1=x_1$, then $y_2=x_2$ and hence Theorem {\ref{s3t1}} reduces to the following result  for the convergence analysis of Iterative Algorithm 2.3 for SpVIP (1.9)-(1.10):

\begin{corollary}\label{s3t3} For each  $i\in\{1,2\},$ let $C_{i}$ be a nonempty, closed and convex subset of $2$-uniformly convex smooth Banach space $E_{i}$ with constant of smoothness $c_i$. Let $F: C_1\to E_1$ be $\alpha_{1}$-strongly monotone  and $\beta_{1}$-Lipschitz continuous, and
let $f: C_2\to E_2$ be $\sigma_{1}$-strongly monotone  and $\eta_{1}$-Lipschitz continuous. Let $A:E_1 \to E_2$  be bounded linear operator. Suppose $x_1 \in C_1 $ is a solution to SpVIP(1.9)-(1.10) then the sequence $\{x_1^n\}$ generated by Iterative Algorithm 2.3  converges strongly to $x_1$ provided that  the constant $ \lambda>0$  satisfies the condition:

$$\left| \lambda-\frac{\alpha_1}{c_1\beta_1^2}\right| < \frac{\sqrt{\alpha_1^2-c_1\beta_1^2(1-p_1^2)}}{c_1\beta_1^2}$$
$$\alpha_1>\beta_1\sqrt{c_1(1 -p_1^2)};~~p_1=\frac{1-m\theta_{2}}{1+m};~~~m=\rho\|A^+\|\|A\|;$$
$$\theta_{2}=\sqrt{1-2\gamma\sigma_1+c_2\gamma^2\eta_1^2};~~~ \gamma > 0.$$
\end
{corollary}

\begin{remark} The extension of the method presented in this paper to split equilibrium problem {\rm\cite{kaz5}} and split variational inclusion {\rm\cite{mou}} needs further research efforts.
\end{remark}

\end{document}